\documentclass[11pt]{amsart}
\usepackage{amssymb,amsmath,amsthm}

\newtheorem{theorem}{Theorem}
\newtheorem{lemma}[theorem]{Lemma}
\newtheorem{cor}[theorem]{Corollary}
\newtheorem{prop}[theorem]{Proposition}

\theoremstyle{definition}
\newtheorem{defn}[theorem]{Definition}

\newcommand{\HH}{\ensuremath{\mathcal{H}}}
\newcommand{\CC}{\ensuremath{\mathbb{C}}}

\newcommand{\RR}{\ensuremath{\mathbb{R}}}
\newcommand{\HL}{\ensuremath{\mathcal{H}L^{2}(\mathbb{C},e^{-\varphi})}}

\title[Pointwise bound for a holomorphic function]
{A POINTWISE BOUND FOR A HOLOMORPHIC FUNCTION
     WHICH IS SQUARE-INTEGRABLE WITH RESPECT TO 
     AN EXPONENTIAL DENSITY FUNCTION}
\author{Kamthorn Chailuek \\
        Wicharn Lewkeeratiyutkul \\}
\address{Department of Mathematics \\
         Faculty of Science \\
         Chulalongkorn University \\
         Bangkok, Thailand 10330}
\email{ckamthorn@hotmail.com, Wicharn.L@chula.ac.th}

\begin{document}

\begin{abstract}
Let $\varphi$ be a real-valued smooth function on $\CC$ satisfying
$0\leq\Delta\varphi\leq M$ for some $M \ge 0$. Denote by $\HL$ the
space of all holomorphic functions which are square-integrable
with respect to the measure $e^{-\varphi(z)}\,dz$. In this paper,
we obtain a pointwise bound for any function in this space. We
show that there exists a constant $K$ depending only on $M$
such that
\[
|f(z)|^2\leq Ke^{\varphi(z)}\|f\|^2_{L^2(\CC,e^{-\varphi})}
\]
for any $f\in\HL$ and any $z\in\CC$.
\end{abstract}

\maketitle
\numberwithin{equation}{section}
\setcounter{page}{1}
\pagenumbering{arabic}

\section{Introduction}

Let $U$ be a non-empty open subset of $\CC$. Denote by $\HH L^2(U,\alpha)$
the space of all holomorphic functions on $U$ which are square-integrable
with respect to the measure $\alpha(\omega)\,d\omega$.

For any $t>0$, consider the Gaussian measure
\[
d\mu_t(z)\;=\;\frac{1}{\pi t}e^{-|z|^2/t}\,dz.
\]
Then the space $\HH L^2(\CC,\mu_t)$ is called the \emph{Segal-Bargmann space}.
See [GM], [H1], [H2], [F] for detailed discussion about the importance of this space,
and its relevance in quantum theory. It is well-known that a pointwise bound
for any function $f \in \HH L^2(\CC,\mu_t)$ is given by
\begin{equation}\label{pwb:s:intro}
   |f(z)|^2 \; \leq \; e^{|z|^2/t}\|f\|^2_{L^2(\CC,\mu_t)}.
\end{equation}
This pointwise bound first appeared in Bargmann's paper [B] and was revisited many times
by other authors. More generally, for any space $\HH L^2(U,\alpha)$, there exists a function
$K(z,\omega)$ on $U\times U$, called the \emph{reproducing kernel}, such that
\begin{equation}\label{pwb:k:intro}
   |f(z)|^2 \; \leq \; K(z,z)\|f\|^2_{L^2(U,\alpha)}
\end{equation}
for any $f\in \HH L^2(U,\alpha)$ and $z\in U$. The Bargmann's pointwise bound
(\ref{pwb:s:intro}) for $\HH L^2(\CC,\mu_t)$ follows from the following formula
of the reproducing kernel for the Segal-Bargmann space:
\begin{equation}\label{rep:sb:intro}
   K(z,\omega) \; = \; e^{z\overline{\omega}/t}.
\end{equation}

In this work, we study a pointwise bound for a function in a more general holomorphic
function space. First, we look at the space $\HH L^2(\CC,e^{-\varphi})$, where
$\Delta\varphi$ is a positive constant. Note that $\Delta(|z|^2/t) = 4/t > 0$,
so this is a generalization of the standard Segal-Bargmann space $\HH L^2(\CC,\mu_t)$.
The technique used here will be that of holomorphic equivalence [H1]. Two holomorphic function
spaces $\HH L^2(U,\alpha)$ and $\HH L^2(U,\beta)$ are holomorphically equivalent if there
exists a nowhere-zero holomorphic function $\phi$ on $U$ such that
\[
   \beta(z)=\frac{\alpha(z)}{|\phi(z)|^{2}} \qquad \text{for all $z\in U$}.
\]
If $\HH L^2(U,\alpha)$ and $\HH L^2(U,\beta)$ are holomorphically equivalent spaces, then
their reproducing kernels are related by
\begin{equation}\label{h:equi:intro}
\alpha(z)K_{\alpha}(z,z)\;=\;\beta(z)K_{\beta}(z,z).
\end{equation}
We show that if $\Delta\varphi=c>0$, then $\HH L^2(\CC,e^{-\varphi})$ is holomorphically
equivalent to the Segal-Barmann space $\HH L^2(\CC,\mu_t)$ where $t=4/c$. It follows
from (\ref{pwb:k:intro}) and (\ref{h:equi:intro}) that
\[
   |f(z)|^2 \; \leq \; \frac{c}{4\pi}e^{\varphi(z)}\|f\|_{L^2(\CC,e^{-\varphi})}^2,
\]
for any $f \in \HH L^2(\CC,e^{-\varphi})$ and any $z\in\CC$.\\
\indent Next, we turn to the space $\HH L^2(\CC,e^{-\varphi})$,
where $\Delta\varphi$ is positive and bounded, i.e.
$0\leq\Delta\varphi\leq M$ for some $M \ge 0$. This space is not
holomorphically equivalent to a Segal-Bargmann space, so we cannot
apply the same technique here.
Our proof relies on a technical lemma which can be stated as follows:
For any $f\in\HH L^2(\CC,e^{-\varphi})$,
\[
|f(0)|^2\;\leq\;Ce^{\varphi(0)}\int_{D(0,1)}|f(\omega)|^2e^{-\varphi(\omega)}\,d\omega
\]
for some $C$ depending only on $M$. By translation to any point $z\in \CC$, we obtain
the following pointwise bound:
\[
   |f(z)|^2 \; \leq \; Ce^{\varphi(z)}\|f\|^2_{L^2(\CC,e^{-\varphi})}.
\]
\indent Here is a brief summary of this work. In section 2, we
study basic properties of holomorphic function spaces. We introduce the concept of
holomorphic equivalence and establish a necessary and sufficient condition for two spaces
to be holomorphically equivalent. In section 3, we establish a pointwise bound for functions
in $\HH L^2(\CC,e^{-\varphi})$.



\section{Holomorphic function spaces}

In this section, we review and prove some relevant facts about holomorphic function spaces
that are needed in this paper. The main reference here is [H1].

\indent Let $U$ be a non-empty open subset of $\CC$. Denote by
$\HH(U)$ the space of all holomorphic functions on $U$. If
$\alpha$ is a strictly positive function on $U$, let
$L^2(U,\alpha)$ be the space of all functions on $U$ which are
square-integrable with respect to the measure $\alpha(\omega)\,d\omega$.
Then $L^2(U,\alpha)$ is a Hilbert space. Let $\HH L^2(U,\alpha)=\HH(U)\cap L^2(U,\alpha)$.
Then $\HH L^2(U,\alpha)$ is a closed subspace of $L^2(U,\alpha)$ and hence a Hilbert space.
Moreover, it is well-known that $\HH L^2(U,\alpha)$ is separable.
\begin{defn}
A \emph{Segal-Bargmann space} is a space $\HH L^2(\CC,\mu_t)$, where
\[
\mu_t(z)=\frac{1}{\pi t}e^{-|z|^2/t}
\]
for some $t>0$.
\end{defn}

Let $K \colon U \times U \to \CC$ be a reproducing kernel for the space $\HH L^{2}(U,\alpha)$.
We refer to [H1] for details of the discussion below. If $\{e_{i}\}_{i=0}^{\infty}$ is
an orthonormal basis for $\HH L^{2}(U,\alpha)$, then the reproducing kernel $K$ is given by
\begin{equation}\label{ker:ort:eqt}
   K(z,\omega) \; = \; \sum_{i=0}^{\infty}e_{i}(z) \overline{e_{i}(\omega)} \qquad (z,\,\omega\in U).
\end{equation}
If we know the reproducing kernel of the space, the pointwise bound of any function $f$ in
$\HH L^{2}(U,\alpha)$ can be obtained by
\begin{equation}\label{est:k:f}
|f(z)|^{2} \; \leq \; K(z,z)\|f\|_{L^2(U,\alpha)} ^{2}.
\end{equation}
Moreover, for a fixed value of $z$, $K(z,z)$ is the smallest constant which makes the pointwise
bound (\ref{est:k:f}) holds for all $f \in \HH L^{2}(U,\alpha)$.

\begin{defn} \label{defn:holo-equiv}
Holomorphic function spaces $\HH L^{2}(U,\alpha)$ and  $\HH L^{2}(U,\beta)$ are
said to be \textit{holomorphically equivalent} spaces if there
exists a nowhere zero holomorphic function $\phi$ on $U$ such that
\[
   \beta(z) = \frac{\alpha(z)}{|\phi(z)|^{2}} \qquad \text{for all $z\in U$.}
\]
In this case, the map $f \mapsto \phi f$ is a unitary map from $\HH L^{2}(U,\alpha)$
onto $\HH L^2(U,\beta)$.
\end{defn}

\begin{lemma}\label{inv:k}
Let \,$\HH L^{2}(U,\alpha)$ and  $\HH L^{2}(U,\beta)$ be holomorphically equivalent
spaces. Let $K_{\alpha}$ and $K_{\beta}$ be their respective
reproducing kernels. Then for each $z\in U$,
\[
\alpha(z)K_{\alpha}(z,z)=\beta(z)K_{\beta}(z,z).
\]
\end{lemma}
\begin{proof}
By formula \ref{ker:ort:eqt} and the fact that a unitary map preserves orthonormal bases,
we obtain
\[
   K_{\beta}(z,\omega) \; = \; \phi(z)\overline{\phi(\omega)}K_{\alpha}(z,\omega).
\]
It follows that
\[
K_{\beta}(z,z)\; = \; |\phi(z)|^{2}K_{\alpha}(z,z) = \;\frac{\alpha(z)}{\beta(z)}K_{\alpha}(z,z).
\]
Thus, \,$\alpha(z)K_{\alpha}(z,z)=\beta(z)K_{\beta}(z,z)$.
\end{proof}
\indent The next goal in this section is to establish a necessary
and sufficient condition for two spaces to be holomorphically
equivalent.

\begin{lemma} \label{Lemma-before-nec-suff}
Let $U$ be an open simply connected set in $\CC$ and $\alpha$ a
strictly positive smooth function on $U$. Then there exists a
holomorphic function $\phi$ such that $|\phi|^{2}=\alpha$ if and
only if \,$\log \alpha$ is harmonic.
\end{lemma}
\begin{proof}
$(\Rightarrow)$ Since $\phi \in \HH(U)$, by a standard result in complex analysis,
there exists a function $\theta \in \HH(U)$ such that $\phi = e^{\theta}$.
Let $u=\text{Re}\,\theta$. Thus, $|\phi|=e^{u}$ and hence $\alpha = e^{2u}$. Then
$\log\alpha = 2u$, which implies that $\Delta \log \alpha = \Delta 2u = 0$.\\
$(\Leftarrow)$ Assume that $u = \log \alpha$ is harmonic. Then there exists
a holomorphic function $f$ such that $u=\text{Re}f$. Hence, $e^{f}$ is also
holomorphic. Let $\phi = e^{f/2}$. Then $\phi \in \HH(U)$ and $e^f = \phi^2$.
Hence, $\alpha = e^{u}=|e^{f}|=|\phi|^{2}.$
\end{proof}

\begin{prop} \label{nec&suf}
Let $U$ be an open simply connected set in $\CC$ and $\alpha$,
$\beta$ strictly positive smooth functions on $U$. Then $\HH
L^{2}(U,\alpha)$ and $\HH L^{2}(U,\beta)$ are holomorphically
equivalent spaces if and only if \,$\Delta \log \alpha(z)=\Delta
\log \beta(z)$.
\end{prop}
\begin{proof}
If $\HH L^{2}(U,\alpha)$ and $\HH L^{2}(U,\beta)$ are holomorphically equivalent, then
there is a function $\phi \in \HH(U)$ such that $\phi \neq 0$ and $|\phi(z)|^{2}=
\frac{\alpha(z)}{\beta(z)}$. By Lemma \ref{Lemma-before-nec-suff},
$\log \frac{\alpha(z)}{\beta(z)}$ is harmonic. Hence, $\Delta(\log\alpha(z)-\log\beta(z))=0$,
which shows that $\Delta\log\alpha(z)=\Delta\log\beta(z)$. It is easy to see that the
reverse implication is true in each step.
\end{proof}
This immediately implies the following corollary:
\begin{cor}\label{cor:nec&suf}
A holomorphic function space $\HH L^{2}(\CC,\alpha)$, where $\alpha$ is a strictly positive
smooth function on $\CC$, is holomorphically equivalent to one of the Segal-Bargmann spaces
if and only if $\Delta \log \alpha =c < 0$. In particular, if $\varphi$ is a smooth function
and $\Delta \varphi $ is a positive constant, then the space $\HH L^2(\CC,e^{-\varphi})$ is
holomorphically equivalent to a Segal-Bargmann space.
\end{cor}
\begin{proof}
Note that if
\[
   \mu_t(z)\;=\;\frac{1}{\pi t}e^{-|z|^2/t},
\]
then
\[
   \Delta\log\mu_t(z)\,=\,-\Delta\frac{|z|^2}{t}
                     \,=\,-\frac{4}{t}\frac{\partial^2}{\partial z \partial\overline{z}}(z\overline{z})
                     \,=\,-\frac{4}{t}\,<\,0.
\]
Thus if \,$\HH L^2(\CC,\alpha)$ is holomorphically equivalent to the Segal-Bargmann space
$\HH L^2(\CC,\mu_t)$, then $\Delta\log\alpha = \Delta\log\mu_t < 0$.

Conversely, if $\Delta\log\alpha = c < 0$, then $\Delta\log\alpha = \Delta\log\mu_t$
where $t=-4/c$. Therefore, $\HH L^2(\CC,\alpha)$ is holomorphically equivalent to the
Segal-Bargmann space $\HH L^2(\CC,\mu_t)$, where $t=-4/c$.
\end{proof}

\section{A pointwise bound for a function in $\HH L^2(\CC,e^{-\varphi})$}

In this section, we obtain a pointwise bound for any function in the holomorphic function space
$\HH L^2(\CC,e^{-\varphi})$. First, we look at the case where $\Delta\varphi$ is a positive
constant.


\begin{theorem}
Let $\varphi$ be a smooth function such that $\Delta\varphi=c$
where $c$ is a positive constant. Then, for any $f\in\HH
L^2(\CC,e^{-\varphi})$ and any $z\in\CC$,
\begin{equation} \label{pointwise-bound-const-case}
   |f(z)|^2 \; \leq \; \frac{c}{4\pi}e^{\varphi(z)}\|f\|_{L^2(\CC,e^{-\varphi})}^2.
\end{equation}
\end{theorem}
\begin{proof}
By Corollary \ref{cor:nec&suf}, $\HH L^2(\CC,e^{-\varphi})$ is
holomorphically equivalent to $\HH L^2(\CC,\mu_t)$, where $t=4/c$.
Then, by Lemma \ref{inv:k},
\[
  K_{e^{-\varphi}}(z,z) \;=\; \frac{1}{\pi t}e^{\varphi(z)} \;=\; \frac{c}{4\pi}e^{\varphi(z)}.
\]
It follows that
\[
   |f(z)|^2 \; \leq \; \frac{c}{4\pi}e^{\varphi(z)}\|f\|_{L^2(\CC,e^{-\varphi})}^2,
\]
for any $f\in\HH L^2(\CC,e^{-\varphi})$ and any $z\in\CC$.
\end{proof}

Note that when $\varphi = |z|^2/t$, we have $c = \Delta \varphi = 4/t$. Hence, in this case
(\ref{pointwise-bound-const-case}) reduces to the usual pointwise bound (\ref{pwb:s:intro})
for the Segal-Bargmann space.


Next, we turn to the situation in which $0 \leq \Delta\varphi \leq M$. The main result is
contained in Theorem \ref{pwb:general-case}. But first we need to establish a technical lemma.

Recall that the function $\Gamma$ defined by
\[
   \Gamma(z)=\frac{1}{2\pi}\log|z|
\]
is the \emph{fundamental solution} for the Laplace's equation on $\RR^2$.
Thus if $\psi\in C_c^{\infty}(\CC)$, then
\[
   \Phi(z) \; = \; \Gamma\ast\psi(z) \; = \; \int_{\CC}\Gamma(\zeta)\psi(z-\zeta)\,d\zeta
\]
satisfies $\Delta\Phi=\psi$.

\begin{lemma}\label{z:0}
Let $\varphi\in C^{\infty}(\CC)$ satisfying $0 \leq \Delta\varphi \leq M$. Then there exists
a constant $C$ depending only on $M$ such that for any $f \in \HH L^2(\CC,e^{-\varphi})$,
\[
   |f(0)|^2 \leq
            C e^{\varphi(0)}\int_{D(0,1)}|f(\omega)|^2e^{-\varphi(\omega)}\,d\omega.
\]
\end{lemma}
\begin{proof}
Choose a function $g\in C_{c}^{\infty}(\CC)$ such that $0 \le g \le 1$,
$g=1$ on $\overline{D(0,1)}$ and $g=0$ outside $D(0,2)$. Let $\psi = g\,\Delta\varphi$.
Then $\psi\in C_{c}^{\infty}(\CC)$, $0 \le \psi \le M$, $\psi=\Delta\varphi$ on
$\overline{D(0,1)}$ and $\psi=0$ outside $D(0,2)$. Thus $\Phi=\Gamma\ast\psi$ satisfies
\begin{equation}\label{con:h:e}
\Delta\Phi(z)=\psi(z)=\Delta\varphi(z)
\end{equation}
for all $z\in D(0,1)$. First, we show that $\Phi$ is bounded above on $D(0,1)$.
Note that $\Gamma(\zeta)\le 0$ if and only if $\zeta \in D(0,1)$.
For any $\omega \in D(0,1)$, we have
\begin{align*}
\Phi(\omega)\;
&= \; \int_{\CC}\Gamma(\zeta)\psi(\omega-\zeta)\,d\zeta \\
&= \; \int_{D(\omega,2)}\Gamma(\zeta)\psi(\omega-\zeta)\,d\zeta \\
&\leq \; \int_{D(\omega,2)\setminus D(0,1)}\Gamma(\zeta)\psi(\omega-\zeta)\,d\zeta \\
&\le \; \frac{M}{2\pi}\int_{D(\omega,2)\setminus D(0,1)}\log|\zeta|d\zeta.
\end{align*}
This shows that $\Phi(\omega) \le BM$ for all $\omega \in D(0,1)$, where
\[
B = \frac{1}{2\pi}\sup_{\omega\in D(0,1)}\int_{D(\omega,2)\setminus D(0,1)}\log|\zeta|d\zeta.
\]
Write $\mathcal{U}= D(0,1)$ and let $h\in\HH L^2(\mathcal{U},e^{-\Phi})$. Fix $0 < s < 1$.
It is not hard to show that
\[
   h(0) = \frac{1}{\pi s^2}\int_{D(0,s)}h(\omega)\,d\omega.
\]
By the Cauchy-Schwarz inequality, it follows that
\[
   |h(0)|^2 \; \leq \; (\pi s^2)^{-2} \left\|\chi_{D(0,s)}e^{\Phi}\right\|_{L^2(\mathcal{U},e^{-\Phi})}^2
                       \|h\|_{L^2(\mathcal{U},e^{-\Phi})}^2.
\]
Hence,
\begin{align*}
   \left\|\chi_{D(0,s)}e^{\Phi}\right\|_{L^2(\mathcal{U},e^{-\Phi})}^2\;
   = \; \int_{D(0,s)}e^{\Phi(\omega)}\,d\omega \;
   \leq \; \int_{D(0,s)}e^{BM}\,d\omega \;
   = \; e^{BM}\pi s^2.
\end{align*}
Thus, for any $0 < s < 1$,
\[
   |h(0)|^2 \; \leq \; \frac{e^{BM}}{\pi s^2}\,\|h\|^2_{L^2(\mathcal{U},e^{-\Phi})}.
\]
It follows that
\[
   |h(0)|^2 \; \leq \; \frac{e^{BM}}{\pi}\,\|h\|^2_{L^2(\mathcal{U},e^{-\Phi})}
\]
for all $h \in \HH L^2(\mathcal{U},e^{-\Phi})$. By a property of the reproducing kernel
(see the paragraph preceding Definition \ref{defn:holo-equiv}) we then have
\[
   K_{e^{-\Phi}}(0,0) \; \leq \; \frac{e^{BM}}{\pi}
\]
where $K_{e^{-\Phi}}$ is the reproducing kernel for $\HH L^2(\mathcal{U},e^{-\Phi})$.

Let $K_{e^{-\varphi}}$ be the reproducing kernel for $\HH L^2(\mathcal{U},e^{-\varphi})$.
Then, by equation (\ref{con:h:e}) and Proposition \ref{nec&suf}, $\HH L^2(\mathcal{U},e^{-\varphi})$
and $\HH L^2(\mathcal{U},e^{-\Phi})$ are holomorphically equivalent and hence,
by Lemma \ref{inv:k},
\begin{align*}
K_{e^{-\varphi}}(0,0) \; = \; \frac{e^{-\Phi(0)}}{e^{-\varphi(0)}}K_{e^{-\Phi}}(0,0) \;
\leq \; C\,e^{\varphi(0)},
\end{align*}
where $C=e^{BM-\Phi(0)}/\pi$.
Thus
\[
   |h(0)|^2 \; \leq \; Ce^{\varphi(0)}\|h\|_{L^2(\mathcal{U},e^{-\varphi})}^2,
\]
for any $h\in\HH L^2(\mathcal{U},e^{-\varphi})$. Let $f \in \HH L^2(\CC,e^{-\varphi})$ and
$h=f\big|_{\mathcal{U}}$. Then $h\in\HH L^2(\mathcal{U},e^{-\varphi})$ and
\begin{align*}
   |f(0)|^2\;&= \; |h(0)|^2 \\
             &\leq \; Ce^{\varphi(0)}\int_{D(0,1)}|h(\omega)|^2e^{-\varphi(\omega)}\,d\omega \\
             &= \; Ce^{\varphi(0)}\int_{D(0,1)}|f(\omega)|^2e^{-\varphi(\omega)}\,d\omega.
\end{align*}
Finally, it remains to show that we can choose a constant $C$ to depend only on $M$.
By straightforward calculations, we have
\[
   \int_{D(0,1)}\Gamma(\zeta)\,d\zeta \; = \; -\frac{1}{4}.
\]
Now, consider
\begin{align*}
\Phi(0)\; = \; \int_{\CC}\Gamma(\zeta)\psi(-\zeta)\,d\zeta \;
          \ge \;\int_{D(0,1)}\Gamma(\zeta)\psi(-\zeta)\,d\zeta \;
          \ge \; -\frac{M}{4}.
\end{align*}
Thus $e^{-\Phi(0)} \le e^{\frac{M}{4}}$, which shows that $C \le \frac{1}{\pi}e^{(B+\frac{1}{4})M}$.
\end{proof}

\begin{theorem} \label{pwb:general-case}
Let $\varphi\in C^{\infty}(\CC)$ with $0\leq\Delta\varphi\leq M$.
Then there exists a constant $C$ depending only on $M$ such that for any
$f \in \HH L^2(\CC,e^{-\varphi})$ and any $z \in \CC$,
\[
   |f(z)|^2 \; \leq \; Ce^{\varphi(z)}\|f\|_{L^2(\CC,e^{-\varphi})}^2.
\]
\end{theorem}
\begin{proof}
Let $z\in\CC$ and $g_z(\omega)=z+\omega$. Then $0\leq\Delta(\varphi\circ g_z)\leq M$.
Let $f \in \HH L^2(\CC,e^{-\varphi})$ and $h = f\circ g_{z}$. Then
$h \in \HH L^2(\CC,e^{-\varphi\circ g_z})$ and by Lemma \ref{z:0},
\begin{align*}
|f(z)|^2 \;
    &= \; |f\circ g_z(0)|^2 \; = \; |h(0)|^2 \\
    &\leq\;Ce^{\varphi \circ g_z(0)}\int_{D(0,1)}|h(\omega)|^2e^{-\varphi\circ g_z(\omega)}\,d\omega\\
    &=\;Ce^{\varphi(z)}\int_{D(0,1)}|f\circ g_z(\omega)|^2e^{-\varphi\circ g_z(\omega)}\,d\omega \\
    &=\;Ce^{\varphi(z)}\int_{D(0,1)}|f(z+\omega)|^2e^{-\varphi(z+\omega)}\,d\omega\\
    &\leq\;Ce^{\varphi(z)}\int_{\CC}|f(\omega)|^2e^{-\varphi(\omega)}\,d\omega \\
    &=\;Ce^{\varphi(z)}\|f\|_{L^2(\CC,e^{-\varphi})}^2.
\end{align*}
\end{proof}

\section*{Acknowledgments}
\noindent
The authors are grateful to Brian Hall for helpful suggestions throughout the process of
this work. We also thank Leonard Gross for useful comments.


\begin{thebibliography}{999}

\bibitem[B]{B}
    V. Bargmann,
    On a Hilbert space of analytic functions and an associated integral transform, Part I,
    \emph{Comm. Pure Appl. Math.} \textbf{14} (1961), 187--214.

\bibitem[F]{F}
    G. Folland, ``Harmonic analysis on phase space,''
    Princeton Univ. Press, Princeton, N.J., 1989.

\bibitem[GM]{GM}
    L. Gross and P. Malliavin,
    Hall's transform and the Segal-Bargmann map, in
    ``Ito's Stochastic Calculus and Probability Theory''
    (M. Fukushima, N. Ikeda, H. Kunita and S. Watanabe, Eds.),
    pp. 73--116. Springer-Verlag, Berlin/New York, 1996.

\bibitem[H1]{H1}
    B. Hall,
    Holomorphic methods in analysis and mathematical physics, in
    ``First Summer School in Analysis and Mathematical Physics''
    (S. P\`{e}rez Esteva and C. Villegas Blas, Eds.), pp. 1--59,
    Contemp. Math., Vol. 260, Amer. Math. Soc., Providence, RI, 2000.

\bibitem[H2]{H2}
    B. Hall,
    Harmonic Analysis with respect to heat kernel measure,
    \emph{Bull. Amer. Math. Soc.}
    \textbf{38} (2001), 43--78.


\end{thebibliography}
\end{document}